\documentclass{svproc}
\usepackage{url}

\usepackage{color}
\usepackage{textcomp}

\usepackage{graphicx}
\usepackage{amssymb}

\begin{document}
\mainmatter              
\title{Algorithmic System Design of Thermofluid Systems}
\titlerunning{Algorithmic System Design of Thermofluid Systems}  
%
\author{Jonas B. Weber \and Ulf Lorenz}
\authorrunning{J.B. Weber \and U. Lorenz} 
%
\tocauthor{Jonas B. Weber, Ulf Lorenz}
\institute{Chair of Technology Management, University of Siegen, \\Unteres Schloss 3, 57072 Siegen, Germany\\
\email{\{jonas.weber, ulf.lorenz\}@uni-siegen.de}}

\maketitle              

\begin{abstract}

Technical components are usually well optimized. However, simply combining these optimized components in a technical system does not necessarily lead to optimal systems. Therefore, focusing on a system perspective reveals new potential for optimization. In this context, we examine thermofluid systems which can be interpreted as fluid systems with superimposed heat transfer. The structure of such systems can be abstracted as a graph -- more specifically, a flow network. We translate the underlying optimization problem into a mixed-integer linear program which is designed to obey the physical laws of heat transfer. Typically, fluid systems can be considered as quasi-stationary systems since their dynamic effects are usually negligible. However, for thermofluid systems this assumption does not hold because time-dependency is an issue as storage tanks for heated fluid gain importance. In order to handle the dynamic effects induced by the storage tanks, we further introduce a continuous-time representation based on a global event-based formulation.
\end{abstract}

\section{Introduction}

Manufacturers of technical components typically optimize their products with regard to a certain operating point. Nevertheless, simply combining these optimized components does not guarantee an optimal system. The design process of technical systems is much more challenging. It involves combining the intended functionality, layout, used components as well as the expected loads for the future use. Empirical studies \cite{VDI} show that these initial decisions make up 70 - 85 \% of the system's total lifespan costs. Therefore, paying attention to a holistic design process reveals new potential. To make use of this potential, we investigate an algorithmic system design approach.

Originated from the German Research Foundation (DFG) founded Collaborative Research Center (SFB) 805 `Control of Uncertainties in Load-Carrying Structures in Mechanical Engineering' at the Technical University of Darmstadt there has been a series of research on the optimization of technical systems over the past years. In this context, mathematical models as well as algorithms for the design and operation of technical systems have been developed. The most extensively investigated topic in this regard and the foundation of this work are flow networks -- more specifically, fluid systems. Therefore, we provide a short overview of the topic. 

In \cite{Martin1} the mathematical optimization of water supply networks has been investigated. In this regard, Morsi et al. \cite{Morsi1} introduced a mixed-integer linear modelling approach based on the piecewise linearization of non-linear constraints for the optimization of dynamic water supply systems with a given layout. Gei{\ss}ler et al. \cite{Geissler1} used a similar approach for the optimization of dynamic transport networks which, in addition to water supply network optimization, has also been applied on the example of transient gas optimization. The investigation of gas networks is carried out further at the SFB-Transregio 154 `Mathematical Modelling, Simulation and Optimization Using the Example of Gas Networks' at the Friedrich-Alexander-University Erlangen-N{\"u}rnberg. Besides optimizing fluid systems with a given layout, F{\"u}genschuh et al. \cite{Fuegenschuh1} examined the optimal layout for the application example of sticky separation in waste paper processing. In this research a mixed-integer non-linear program (MINLP) for the simultaneous selection of the network topology as well as the optimal settings of each separator for the steady state has been proposed. A MINLP formulation to design decentralized water supply systems for skyscrapers has been used in \cite{Leise}.  Furthermore, P{\"o}ttgen et al. \cite{Poettgen2} compared linear and non-linear programming techniques for the combined layout and control optimization of booster stations, while favorable combinations between model formulations and mathematical solver packages have been studied in \cite{Rausch}. In addition to standard solver packages, problem specific primal and dual solution algorithms for the linear formulation have been developed in \cite{Weber} to speed up the solution process. 

At the same time, a separate research area arose. This research area, called `Technical Operations Research' (TOR), combines technical and mathematical know-how in order to design optimal technical systems. In terms of content and objective, it is close to current research. However, the focus is on providing engineers with tools that enable the use of modern mathematical methods during the system design process in the form of applicable software. In this context, Pelz et al. \cite{Pelz1} introduced an `artificial fluid system designer' as an attempt to automatically find optimal pump system designs. Beyond fluid systems, this approach has also been applied to the optimization of other technical systems such as gearboxes \cite{Altherr2} and lattice structures \cite{Reintjes}.

In addition to  the transport of fluid, many technical applications involve heating and cooling. The corresponding systems are then called thermofluid systems. Though these systems can be regarded as fluid systems with superimposed heat transfer from a technical point of view, both parts have to be considered at the same time since they mutually depend on each other. In this spirit, an already existing heating circuit of a conference center in Darmstadt has been examined and optimized in \cite{Poettgen4}. In this paper, however, we present a mixed-integer linear program (MILP) to model thermofluid systems based on a graph representation for system synthesis tasks. Besides that, another aspect of heat transfer is considered. While it can be reasonable for fluid systems to be assumed as quasi-stationary systems, this does not hold for most thermofluid systems. Although some can be regarded as quasi-stationary, many of them show dynamic behavior. A major reason for this is that when dealing with heat, storage tanks become important. Thus, an appropriate time representation in order to model the dynamic behavior is required. For this purpose, we introduce a continuous-time approach.

\section{Modeling of Thermofluid Systems}


The system synthesis task considered in this paper can be stated as follows: Given a construction kit of technical components such as pumps, pipes or boilers as well as a technical specification of load collectives, compare all valid systems and choose the one for which the lifespan costs  -- the sum of purchase costs and the expected energy costs -- are minimal. In this context, a system is called a valid system if it is able to satisfy every prospected load.

A possible representation for this purpose is to model a thermofluid system as a source-target-network $(G,s,t)$ with directed multigraph $G := (V,E)$. $E$ is the set of edges representing technical components. The set of vertices $V$ represents interconnections between components, whereas $s,t \in V$ are two distinguished vertices, namely the source and the target of the network. Each possible system can be represented as a connected subgraph without directed cycles defined by a purchase decision of components -- indicated by binary decision variables. 

\subsection{Basic Network Model for Fluid Systems}

\begin{table}[!htb]
\caption{Variables and Parameters}
\label{vars_and_pars}
\begin{center}
\begin{tabular}{c@{\quad}cl}
\hline
\multicolumn{1}{c}{\rule{0pt}{12pt}
                   Symbol}&\multicolumn{1}{c}{\rule{0pt}{12pt}Range}&\multicolumn{1}{c}{\rule{0pt}{12pt}Description}\\[2pt]
\hline\rule{0pt}{12pt}
$b_{i,j}$ & $\{0;1\}$ & Purchase decision of component $(i,j)$\\
$a_{i,j}^s$ & $\{0;1\}$ & Activation decision of component $(i,j)$ in scenario $s$\\
$\dot{v}_{i,j}^s$ & $\mathbb{R}_0^+$ & Volume flow through component $(i,j)$ in scenario $s$ \\
$\dot{q}_{i,j}^{in \; s}$ & $\mathbb{R}_0^+$ & Heat flux at the inlet of component $(i,j)$ in scenario $s$ \\
$\dot{q}_{i,j}^{out \; s}$ & $\mathbb{R}_0^+$& Heat flux at the outlet of component $(i,j)$ in scenario $s$ \\
$h_k^s$ & $\mathbb{R}_0^+$ & Pressure head at connection $k$  in scenario $s$ \\
$t_k^s$ & $\mathbb{R}_0^+$ & Temperature at connection $k$  in scenario $s$ \\
$T_{i,j}(\dot{v}_{i,j}^s,\dot{q}_{i,j}^{in/out \; s})$ & $\mathbb{R}_0^+$ &Inlet/outlet temperature of component $(i,j)$ in scenario $s$ \\
$P_{i,j}(\dot{v}_{i,j}^s,n_{i,j}^{s})$ & $\mathbb{R}_0^+$ & Power consumption of pump $(i,j)$ in scenario $s$ \\
$\Delta \dot{q}_{i,j}^s$ & $\mathbb{R}$ & Heat increase by heat source $(i,j)$ in scenario $s$ \\
$t_{i,j}^s$ & $\mathbb{R}_0^+$ & Outlet temperature of temp. source $(i,j)$ in scenario $s$ \\
$\Delta h_{i,j}^s$ & $\mathbb{R}$ & Pressure increase by component $(i,j)$ in scenario $s$ \\
$ l_{i,j}^s$ &  -  & Level variable of component $(i,j)$ in scenario $s$\\ 
\hline\rule{0pt}{12pt}
$\dot{V}^{max}$ & - & Upper bound on the volume flow\\
$\dot{Q}^{max}$ & - & Upper bound on the heat flux\\
$H^{max}$ & - & Upper bound on the pressure head\\
$T^{max}$ & - & Upper bound on the temperature\\
$\Delta \dot{Q}_{i,j}^{max/min}$ & - & Max./min. heat increase by heat source $(i,j)$\\
$T_{i,j}^{max/min}$ & - & Max./min. outlet temperature of temp. source $(i,j)$\\
$C^{kWh}$ & - & Energy costs per kilowatt hour\\
$C^{buy}_{i,j}$ & - & Purchase costs of component $(i,j)$ \\
$E_{i,j}$ & - & Efficiency of heating or cooling component $(i,j)$\\
$OLS$ & - & Operational lifespan of the system\\
$F^s$ & - & Fraction of the operational lifespan $OLS$ of scenario $s$ \\[2pt]
\hline
\end{tabular}
\end{center}
\end{table}

Since thermofluid systems are an extension of fluid systems, we shortly introduce the relevant constraints for the quasi-stationary case. For a more detailed view of the underlying logical, physical and technical properties, we refer to \cite{Weber}. An overview of all variables and parameters used is given in Table \ref{vars_and_pars}. Note that because of the quasi-stationarity, similar loads of a load collective are aggregated to so-called load scenarios $S$. These load scenarios occur for a specific portion of the total time. Two physical quantities are necessary in order to describe a fluid system, namely the extensive flow variable volume flow $\dot{V}$ on edges and the intensive variable pressure or more explicitly, pressure head $H$ on vertices.

A component can only be used to satisfy a load scenario if it is installed:

\begin{equation}\label{FluidNB1}
\forall \, s \in S, \; \forall \, (i,j) \in E: \quad a_{i,j}^{s} \leq b_{i,j} 
\end{equation}

The pressure head at each connection must be reasonable:

\begin{equation}\label{FluidNB2}
\forall \, s \in S, \; \forall \, k \in V: \quad h_k^s \leq H^{max} 
\end{equation}

If a component is operational, Bernoulli's equation applies. In case of pumping components, e.g. pumps, the pressure increase caused by the component increases the pressure at its outlet and therefore the adjacent system pressure. For non-pumping components the pressure increase $\Delta h_{i,j}^s$ is typically 0:


\begin{eqnarray}\label{FluidNB3}
\forall \, s \in S, \; \forall \, (i,j) \in E: \quad h_j^s - h_i^s \leq \Delta h_{i,j}^s + H^{max} \cdot (1 - a_{i,j}^s) \\
\quad h_j^s - h_i^s \geq \Delta h_{i,j}^s - H^{max} \cdot (1 - a_{i,j}^s)
\end{eqnarray}

If a component is operational, its volume flow must be reasonable, otherwise it vanishes:


\begin{equation}\label{FluidNB4}
\forall \, s \in S, \; \forall \, (i,j) \in E: \quad \dot{v}_{i,j}^s \leq  \dot{V}^{max} \cdot a_{i,j}^s
\end{equation}

For all vertices, except for the source and the sink, the continuity equation applies which demands for flow conservation:

\begin{equation}\label{FluidNB5}
\forall \, s \in S, \; \forall \, k \in V \, \backslash \{ s,t\}: \quad \sum_{(i,k) \in E}\dot{v}_{i,k}^s - \sum_{(k,j) \in E}\dot{v}_{k,j}^s = 0
\end{equation}

In addition, if a component is designated as a pump, the pump's operating point must lie on its characteristic curve. This can be achieved by generating a suitable number of base points from the empirically known $H_{i,j}(\dot{v}_{i,j}^s, n_{i,j}^s)$ and $P_{i,j}(\dot{v}_{i,j}^s, n_{i,j}^s)$ functions and forcing the respective variables on the linearized curves defined by these points. The used linearization techniques follow \cite{Vielma}.

\subsection{Heat Transfer}

In order to describe the properties related to heating and cooling three physical quantities are necessary. These are the already introduced volume flow $\dot{V}$, the heat flux $\dot{Q}$ as another extensive flow variable on edges as well as the intensive variable temperature $T$ on vertices.

All three are coupled by the specific heat formula. In this connection, the specific heat is the amount of heat per unit mass required to raise the temperature by one Kelvin  -- or degree Celsius. The relationship is typically expressed as shown in Equation \ref{specific_heat} where c is the specific heat. As an example, the specific heat of water -- the common substance with the highest specific heat -- is one calorie per gram and degree Celsius which is equal to 4.184 joule per gram and degree Celsius at a temperature of 18\textdegree C. However, the relationship does not hold if phase changes occur due to the fact that heat added or removed during a phase change does not change the temperature. 

\begin{equation}
\label{specific_heat}
\Delta Q = m \cdot c \cdot \Delta T
\end{equation}

Keeping this in mind the equation can be rewritten and simplified in order to obtain the required flow variables -- indicated by the dot notation -- by assuming that the used fluid is water with a density of about 1 and referring all values to a reference temperature of 0\textdegree C (273.15 K):

\begin{equation}
\label{specific_heat_simple}
\dot{Q} = \dot{V} \cdot c \cdot T
\end{equation}

Another important relation exists if we intend to mix (possibly different) fluids with different temperatures. In this case, the mixing temperature $T_{M}$ can be easily calculated using the formula shown in Equation \ref{mixing} where $|N|$ fluids are mixed.
\begin{equation}
\label{mixing}
T_M = \frac{\sum\limits_{i \in N} m_i \cdot c_i \cdot T_i}{\sum\limits_{i \in N} m_i \cdot c_i}
\end{equation}

Once again, we can simplify this relation by assuming that we mix the same kind fluid -- in this case water with a density of 1 -- and rewrite it for flow variables:

\begin{equation}
\label{mixing_simple}
T_M = \frac{\sum\limits_{i \in N} \dot{V}_i \cdot T_i}{\sum\limits_{i \in N} \dot{V}_i} = \frac{\sum\limits_{i \in N} \dot{Q}_i}{\sum\limits_{i \in N} \dot{V}_i} \cdot \frac{1}{c}
\end{equation}

Using these basic principles of heat transfer, we can define the required constraints for the model. However, since the heat flux, unlike the volume flow, can change along edges, two variables rather than one are required to model it. The variable $\dot{q}^{in}$ represents the heat flux directly behind a vertex corresponding to a component's inlet. Following this, $\dot{q}^{out}$ represents the heat flux directly in front of a vertex corresponding to a component's outlet.

If a component is operational, its heat flux must be reasonable, otherwise it vanishes:

\begin{eqnarray}
\label{heat2}
 \forall \, s \in S, \; \forall \, (i,j) \in E: \quad \dot{q}_{i,j}^{in  \; s} \leq \dot{Q}^{max} \cdot a_{i,j}^s \\
 \quad \dot{q}_{i,j}^{out  \; s} \leq \dot{Q}^{max} \cdot a_{i,j}^s
\end{eqnarray}

The heat flux must be preserved at all vertices, except for the source and the sink. This is due to the law of energy conservation:

\begin{equation}
\label{heat1}
\forall \, s \in S, \; \forall \, k  \in V \, \backslash \{ s,t\}: \quad \sum_{(i,k) \in E} \dot{q}_{i,k}^{out \; s} - \sum_{(k,j) \in E} \dot{q}_{k,j}^{in  \; s} =0
\end{equation}

If a heating or cooling component, e.g. boiler, is operational the transferred heat between the component and the fluid affects the heat flux at the component's outlet. For non-heating or -cooling components the difference $\Delta \dot{q}_{i,j}^s$ is typically 0, expect for possible losses:

\begin{eqnarray}
\label{heat3}
\forall \, s \in S, \;  \forall \, (i,j) \in E: \quad \dot{q}_{i,j}^{out  \; s} \leq \dot{q}_{i,j}^{in  \; s} + \Delta \dot{q}_{i,j}^s + (1 - a_{i,j}^s) \cdot \dot{Q}^{max} \\
 \label{heat4}
 \quad \dot{q}_{i,j}^{out  \; s} \geq \dot{q}_{i,j}^{in  \; s} + \Delta \dot{q}_{i,j}^s - (1 - a_{i,j}^s) \cdot \dot{Q}^{max}
\end{eqnarray}



For the purpose of mixing different flows at a vertex, Equation \ref{mixing_simple} has to be met. Because of the bi-linear relationship arising from the specific heat formula, shown in Equation \ref{specific_heat_simple}, the resulting temperature -- depending on the sum of incoming volume flows and heat fluxes -- must be linearized: 

\begin{eqnarray}
 \label{temp2}
\forall \, s \in S, \; \forall \, k \in V: \quad & t_k^s= \frac{\sum\limits_{(i,k) \in E} \dot{q}_{i,k}^{out  \; s}}{\sum\limits_{(i,k) \in E} \dot{v}_{i,k}^s} \cdot \frac{1}{c} \nonumber \\
& = T(\sum\limits_{(i,k) \in E} \dot{v}_{i,k}^s, \sum\limits_{(i,k) \in E} \dot{q}_{i,k}^{out  \; s}) 
\end{eqnarray}

Due to the temperature's property to be an intensive variable -- it does not depend on the amount of the substance for which it is measured -- all flows leaving a vertex must have the same temperature. Again, the bi-linear relationship between the temperature, the outgoing volume flow and the heat flux has to be linearized:

\begin{eqnarray}
 \label{temp3}
\forall \, s \in S, \; \forall \, (i,j) \in E: \quad T(\dot{v}_{i,j}^s, \dot{q}_{i,j}^{in  \; s})  \leq t_i^s + (1 - a_{i,j}^s) \cdot T^{max} \\
\quad T(\dot{v}_{i,j}^s, \dot{q}_{i,j}^{in  \; s})  \geq t_i^s - (1 - a_{i,j}^s) \cdot T^{max}
\end{eqnarray}

\subsection{Heating and Cooling Components}

There are a multitude of different heating and cooling components. Modeling all of them would be time-consuming with a small benefit. However, for the purpose of this paper they can be grouped according to certain characteristics.

A key differentiator in this model is the working principle, i.e. how the heating and cooling tasks are performed. In this context, two ideal sources of thermal energy can be distinguished. Even though real technical components are not ideal thermal sources, most can be assumed as ideal and simplified without much loss. The two ideal sources of thermal energy are: ideal heat sources and ideal temperature sources. An ideal heat source is able to deliver a constant, predefined heat flux  independent of the temperature difference between its inlet and outlet as well as the volume flow. An example for components which are modeled as ideal heat sources in this model are (tankless) boiler. An ideal temperature source in contrast can maintain a predefined temperature at its outlet independent of the heat flux required as well as the inlet temperature and the volume flow. It therefore produces a constant absolute temperature. In the case of ideal temperature sources, Equations \ref{heat3} and \ref{heat4} do not apply. Rather a constant temperature is assigned to the component's outlet as shown in Equations \ref{exchange1} and \ref{exchange2}. Components which can be modeled as an ideal temperature source are heat exchangers for district heat \cite{Poettgen4}:

\begin{eqnarray}
\label{exchange1}
 \forall \, s \in S, \; \forall \, (i,j) \in temp. \, source(E): \quad t_j^s \leq t_{i,j}^s + (1 - a_{i,j}^s) \cdot T^{max}  \\
 \label{exchange2}
 \quad t_j^s \geq t_{i,j}^s - (1 - a_{i,j}^s) \cdot T^{max}
\end{eqnarray}

Another important characteristic is the control of the components. A distinction can be made between components with single-stage, multi-stage or continuously variable control. To model this behavior a maximum heat flux $\Delta \dot{Q}_{i,j}^{max}$ or temperature $T_{i,j}^{max}$, a minimum heat flux $\Delta \dot{Q}_{i,j}^{min}$ or temperature $T_{i,j}^{min}$ and an additional level variable $l_{i,j}^s$ are required as shown in Equation \ref{level_heat} and \ref{level_temp}, respectively. For components with a  continuously variable control $l_{i,j}^s$ is continuously in the range of 0 to 1 whereas $l_{i,j}^s$ is limited to certain fixed values between 0 and 1 for multi-stage controlled components. In the case of single-stage controlled components with only one fixed operating point, the level variable $l_{i,j}^s$ is fixed to 0 or 1:

\begin{equation}
\forall \, s \in S, \; \forall \, (i,j) \in E: \quad a_{i,j}^s \geq l_{i,j}^s
\end{equation}

\begin{equation}\label{level_heat}
\forall \, s \in S, \; \forall \, (i,j) \in E: \quad \Delta \dot{q}_{i,j}^s = \Delta \dot{Q}_{i,j}^{min} + (\Delta \dot{Q}_{i,j}^{max} - \Delta \dot{Q}_{i,j}^{min}) \cdot l_{i,j}^s
\end{equation}

\begin{equation}\label{level_temp}
\forall \, s \in S, \; \forall \, (i,j) \in E: \quad t_{i,j}^s = T_{i,j}^{min} + (T_{i,j}^{max} - T_{i,j}^{min}) \cdot l_{i,j}^s
\end{equation}

Additionally, the energy source used has to be considered as it has an impact on the cost components in the objective function. Typical energy sources used for heating and cooling  are electricity, gas, biomass, (heating) oil, solar energy, district heat or geothermal energy. Due to the conversion of different units of measurement and costs per unit each of them represents a different cost component in the objective function.

Lastly, it is important whether components can be classified as flow or storage components. In general, three groups can be distinguished among the components associated with heating and cooling. Tankless components, components with integrated storage and pure storage components. However, only tankless components are compatible with a quasi-stationary view. If storage components are involved, dynamic systems arise which show time-dependency. A possible handling of such systems is examined in section \ref{dynamic}.

\subsection{Objective}

The objective function contains two parts, the investment costs for the system components and their energy costs. The aim is to minimize the sum of both over an a-priori defined operational lifespan $OLS$. Each load scenario occurs in a fraction $F^s$ of the total time. The investment costs include the purchase costs of the installed components. The energy costs are further subdivided based on the energy source used. For most energy sources the energy terms look alike. Electricity, however, is a special case as it can be used for two purposes, to operate the pumps and thus to convey fluid as well as for heating or cooling tasks. Furthermore, to determine the actual energy consumption, the heat supplied or removed from the system has to be adjusted by the ratio of useful heating or cooling provided to work required, denoted as $E_{i,j}$:

\begin{eqnarray}
\label{objective}
&minimize \quad c_{total}&  \nonumber \\ [1ex]
&c_{total} = c_{invest} + c_{energy}& \nonumber \\ [1ex]
&c_{invest} = \sum (C^{buy}_{i,j} \cdot b_{i,j}) &\nonumber \\ [1ex]
&c_{energy} = c_{electr.} + c_{gas} + c_{biomass} + c_{oil}  + ... + c_{district} &\nonumber \\ [1ex]
&c_{electr.} = C^{kWh}_{electr.} \cdot OLS \cdot \sum\limits_{s \in S} (F^s \cdot (\sum\limits_{(i,j) \in pumps(E)}P_{i,j}(\dot{v}_{i,j}^s, n_{i,j}^s ) &\nonumber  \\ [1.5ex]
& + \sum\limits_{(i,j) \in electr. \, heat(E)}\Delta \dot{q}_{i,j}^s \cdot \frac{1}{E_{i,j}})) & \nonumber \\ [1.5ex]
&...&\nonumber  \\ [1.5ex]
&c_{district} = C^{kWh}_{district} \cdot OLS \cdot \sum\limits_{s \in S} (F^s \cdot (\sum\limits_ {(i,j) \in district(E)}(\Delta \dot{q}_{i,j}^s \cdot \frac{1}{E_{i,j}}))&
\end{eqnarray}


\section{Time Representation for Dynamic Behavior}\label{dynamic}
As indicated earlier, the quasi-stationary view is not applicable whenever similar loads cannot be aggregated to scenarios. This happens whenever the actual state of the system depends on its load history and therefore a path dependency occurs. Regarding thermofluid systems this is the case for systems with storage, existing components with extensive start-up and run-down phases, general delayed system responses or the like. Note, that this section focuses on storage components while it may be adoptable for other purposes, too.

For the representation of time, two contrary types exist -- discrete and continuous representations. The first one divides the observation period into uniform time intervals. All system events -- the internal and external actions that cause the system to leave the stationary state -- are associated with the start or the end of an interval. While the benefits of this representation -- including a reference grid for all operations, an easy implementation and typically well-structured mathematical problems -- seem attractive for some cases, it also has major disadvantages. Because of the a-priori fixed intervals and interval lengths, events are limited to these points in time. For this reason, the discrete representation is only an approximation, with its resolution depending on the number of intervals. However, more intervals lead to higher computational effort. Therefore, a trade-off between accuracy and the computational effort required must be made. Additionally, the discrete representation leads to larger instances than necessary since the intervals must be uniform and therefore the length of an interval is the smallest common divider of the duration each considered (constant) load occurs. This is especially the case with real-world applications.

Due to the discussed disadvantages, we focus on a continuous-time representation. For this, a global event-based approach is used. This means that the event points (or actions) define a joint, unified reference grid for all components of the system while a unit-specific event-based approach would introduce its own reference grid for each component. The basic idea is that (additional) variables are used to determine the timings of the intervals. However, there are also challenges for this approach. Non-linear programs (NLP) arise due to the fact that the interval lengths are no longer constant but variable. Furthermore, the estimation and adjustment of the number of time intervals is a challenge. If the number of intervals is underestimated, inaccurate solutions or even infeasibility may occur. If, on the other hand, the number of intervals is overestimated, unnecessarily big instances arise.

In order to describe the approach, a short introduction to the properties of storage tanks is given at this point. Generally, the filling level of a storage tank at a point in time $t$ can be determined using the flow balance equation:

\begin{equation}
\label{balance}
V_{t} = V_{t-1} + \int_{t-1}^{t} (\dot{V}_{in} - \dot{V}_{out})
\end{equation}

For constant flows between $t-1$ and $t$ the equation simplifies to:

\begin{equation}
V_{t} = V_{t-1} + (\dot{V}_{in} - \dot{V}_{out}) \cdot (\tau_{t} - \tau_{t-1})= V_{t-1} + \Delta \dot{V}  \cdot \Delta \tau
\end{equation}

It can be seen that a non-linear term still exists since the two variables $\Delta \dot{V}$ and  $\Delta \tau$ are multiplied. Nevertheless, the relation becomes easier to handle compared to Equation \ref{balance} if only constant flows occur. The resulting challenge is how to choose the (maximum) number of intervals to ensure that only constant flows occur.

In the following, only one source and one sink are used for illustration purposes. In this case, flows are constant as long as the demand of the system $\dot{V}_{system}$, which corresponds to the demand at the source $\dot{V}_{source}$, is constant:

\begin{equation}
\dot{V}_{system} = \dot{V}_{sink} + \sum_{i \in Tanks} (\dot{V}_i^{in} - \dot{V}_i^{out}) = \dot{V}_{source}
\end{equation}

The demand of the system $\dot{V}_{system}$ always changes when an activity at the sink $\dot{V}_{sink}$, i.e. on the consumer side, takes place or the demand changes indirectly due to the filling or emptying of a storage tank. The change in demand due to the first is called a main-event. While the number of main-events is known in advance because of the a-priori determined projected demands by the consumer -- comparable to the load scenarios in the quasi-stationary case but with respect to their chronological order -- the number of intervals between the main-events still needs to be determined.

If there is a constant demand at the sink, a storage tank should strive to empty as early as possible and to fill as late as possible during this period to avoid energy losses. Even if energy losses are not explicitly considered, it is reasonable to assume that the filling or emptying takes place in only one continuous process instead of multiple, interrupted processes right before or after a main-event without loss of generality. We therefore define that at most one filling and one emptying process per tank takes place between two main events. Thus, the upper bound on the number of intervals between two main-events results from:

\begin{equation}
n_{intervals} = n_{sources} + 2 \cdot n_{tanks}
\end{equation}

\begin{center}
\begin{figure}
\label{events}
\includegraphics[width=\textwidth]{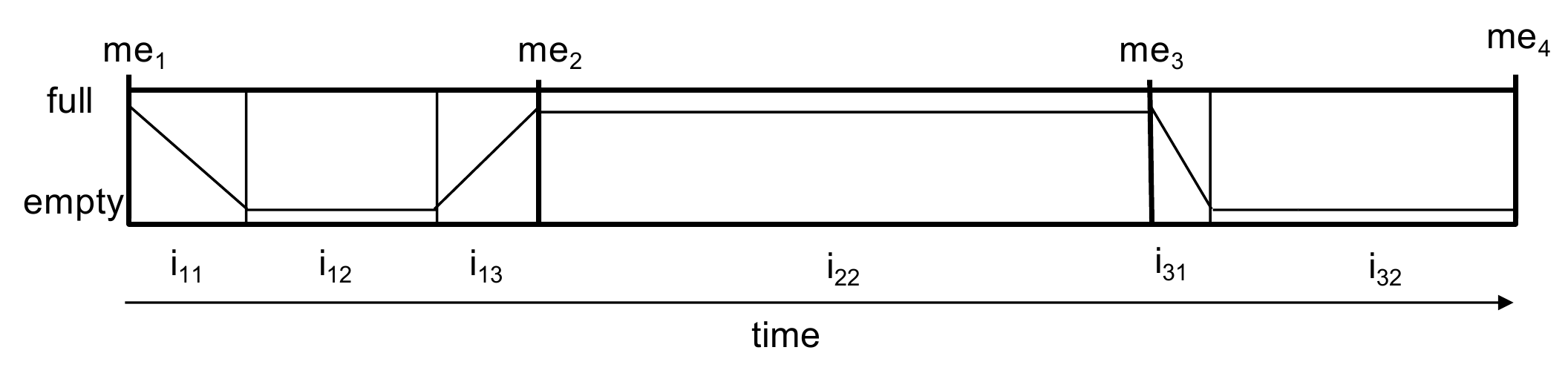}
\caption{Schematic diagram of the filling level of a storage tank over time.}
\end{figure}
\end{center}

In the case of one source and one storage tank, as shown in figure \ref{events}, there are at most three intervals $i$ between two main-events $me$, one for the emptying of the tank $i_{11}$, one if there is no change for the tank $i_{12}$ and one for the filling of the tank $i_{13}$. Figure \ref{events} also illustrates that the determined number of intervals is only an upper bound. Between the main-events $me_1$ and $me_2$ all three intervals are needed, while between $me_2$ and $me_3$ one interval would be sufficient as the tank is neither filled nor emptied and the demand is satisfied by a continuous flow from the source to the sink alone all the time. Therefore, $i_{21}$ and $i_{23}$ exist only theoretically and have a length of 0.

\section{Conclusion and Outlook}

In this paper, we presented an optimization model for the algorithmic system design of thermofluid systems. The model provides a unified framework for later work on technical applications in this field. A possible real-world application are temperature control systems which are used, amongst others, for injection molding. This can be accompanied by an extension to include losses. In addition, we work towards making the model more efficient with regard to the number of linearizations required. One possible measure could be to only allow the mixing of flows for certain connections which would reduce the number of linearizations, though modelling becomes more challenging. In this context, another measure would be to limit the temperature observation to thermal components only. For solving the problem efficiently, we plan to adopt the approach for fluid systems suggested in \cite{Weber} which uses specified primal and dual solution techniques at the same time.

Furthermore, we introduced a continuous time-representation for technical fluid-based systems. The proposed representation enables the consideration of dynamic effects for the optimization of these systems since it is a reasonable compromise between accuracy and computational effort. The dynamic effects are considered in such a way that the essential properties of storage components can still be taken into account while it is simple enough to be applicable for optimization. In a consecutive step, we therefore work on the formulation of an optimization model based on the presented considerations to enable the system synthesis of thermofluid systems with dynamic behavior. The continuous representation is particularly advantageous in comparison to a discrete representation if the number storage tanks is small. In this case, there is only a limited number of linearizations although a high resolution is achieved. This advantage is further enhanced for smaller instances with a manageable number of load changes. A promising application could be to apply the representation to industrial processes with a, e.g. daily, repeating production sequence. Besides that, the upper bound on the number of intervals between two main-events could be further reduced depending on the particular characteristics of the examined application.

%
%

\end{document}